\documentclass[11pt]{article}
\usepackage{latexsym}
\usepackage{amsfonts}
\usepackage{mathrsfs,amsthm,amsmath}
\usepackage{color}
\usepackage{todonotes}
\usepackage{comment}
\newtheorem{theorem}{Theorem}[section]
\newtheorem{proposition}{Proposition}[section]
\newtheorem{definition}{Definition}[section]

\newtheorem{corollary}{Corollary}[section]
\newtheorem{lemma}{Lemma}[section]
\newtheorem{remark}{Remark}[section]
\newtheorem{example}{Example}[section]

\newcommand{\ri}{{\rm i }}
\newcommand{\re}{{\rm e }}
\newcommand{\rd}{{\rm d }}
\newcommand{\B}{\hfill $\Box$}

\begin{document}\title{ On the  connection between the Beurling-Malliavin density and   the asymptotic density}
 \author{Rita Giuliano\thanks{Dipartimento di
		 Matematica, Universit\`a di Pisa, Largo Bruno
		 Pontecorvo 5, I-56127 Pisa, Italy (email: \texttt {rita.giuliano@unipi.it}, ORCID: 0000-0002-1638-2651)}
 \and~Georges Grekos\thanks{Faculté des Sciences et Techniques,  Universit\'e Jean Monnet,  23  rue Dr Paul Michelon. 42023 Saint-Etienne Cedex 2, France (email:\texttt {grekos@univ-st-etienne.fr}, ORCID: 0000-0002-3681-8336).    
 } \footnote {The authors are grateful to the CIRM of Marseille for the hospitality offered to them  in 2019 (under the Program {\it Research in pairs} \texttt {https://conferences.cirm-math.fr/2019-calendar.html}), when this research was at the beginning.}
 }
  
\date{}
\maketitle
\begin{abstract}\noindent We study the notion of Beurling-Malliavin density from the point of view of Number Theory. We prove a general relation between the Beurling-Malliavin density and the upper asymptotic density; we identify a class of sequences for which the two densities coincide; this class contains the arithmetic progressions. Last, by means of  an alternative definition of Beurling-Malliavin density,  we study the connection with the asymptotic density for another kind of sequences that again generalizes the arithmentic progressions.  

 \medskip
\noindent\emph{Keywords}: Beurling-Malliavin density, asymptotic density, increasing sequence, substantial interval, complete system, first digit set, regularly varying function, set of powers

 \medskip
\noindent
\noindent\emph{2020 Mathematical Subject Classification}: 11B25,11B05
%60F10 Large deviations
%60G70 Extreme value theory; extremal stochastic processes
%60F05 Central limit and other weak theorems 
\end{abstract}

 	 ~

 \section{Introduction} Let $ \Lambda= \big(\lambda_n\big)_{n \in \mathbb{Z}}$ be a sequence of real numbers. In order to solve the problem of finding the radius of completeness of $ \Lambda$ (denoted by $\mathcal{R}(\Lambda)$),   A. Beurling and P. Malliavin introduced for the first time in the paper \cite{BM} the quantity $b(\Lambda)$, defined as follows: if $\big(I_n \big)_{n \in \mathbb{Z}}$ is a sequence of disjoint intervals on $\mathbb{R}$,   call it {\it short} if
$$\sum_{n \in \mathbb{Z}}\frac{|I_n|^2}{1 + {\rm dist}^2(0, I_n)}< \infty   $$ 
and {\it long} otherwise; then deﬁne
$$b(\Lambda):=\sup\{d :\exists \,\,{\rm long} \,\,(I_n)_{n\in \mathbb{Z}} \,\,{\rm such\,\, that} \,\, \#(\Lambda \cap I_n) > d|I_n|, \forall n\in \mathbb{Z}\}$$
(the particular formulation used above comes from \cite{P1} and \cite{P2}). The celebrated Theorem of  \cite{BM}  states that the radius of completeness of $ \Lambda$ is connected with $b(\Lambda)$ by the formula
$$\mathcal{R}(\Lambda) = 2 \pi b(\Lambda).$$
Later  the same quantity $b(\Lambda)$ was studied by other authors,  for instance  \cite{Ka},  \cite{K2}, and used  in several papers; we cite   \cite{Kab}, \cite{Kr}  and  \cite{R}; for an extensive study and generalizations see \cite{P}, \cite{P1}, \cite{XY} (see  \cite{P2} for a complete list of references); in particular the author of \cite{R} discovers an alternative equivalent definition, which will be stated and used later (see Definition \ref{def1}).  

\bigskip
\noindent
In the whole literature on the subject the quantity $b(\Lambda)$ is called the \lq\lq density\rq\rq \ of $\Lambda$; in the sequel we call it \lq \lq Beurling-Malliavin density\rq\rq (or \lq \lq BM-density\rq\rq \ for short); anyway, 
to the best of our knowledge, no paper addresses the problem of investigating whether it \lq\lq deserves\rq\rq (so to say) the name of  \lq\lq density\rq\rq , in the sense in which this term is used in Number Theory. Just to give a trivial example, we have not found any paper stating that the BM-density of the set of multiples of the positive integer $p$ is $\frac{1}{p}$.

\bigskip
\noindent
Thus, the most natural notion of density of a sequence $ \Lambda$ of positive integers used in Number Theory  being the  asymptotic density  (denoted     by  the symbol $d(\Lambda)$), the aim
of the present investigation is to study the connection between  $b(\Lambda)$ and $d(\Lambda)$.

 \bigskip
  \noindent
  The Beurling-Malliavin   density is defined in \cite{BM}  for general sequences  $ \Lambda= (\lambda_n)_{n \in \mathbb{Z}}$ of real numbers; on the other hand, if $  \Lambda= (\lambda_n)_{n \in \mathbb{N}^*}$ is a  sequence of integers, we may be interested in  calculating its  asymptotic density. Thus, in order to compare these two concepts,  in what follows we shall confine ourselves to    sequences  of real numbers indexed by $n \in \mathbb{N}^*$, positive, ultimately strictly increasing  and such that  $\lim_{n \to \infty}\lambda_n= + \infty$. In the sequel, by the term {\it sequence} we always mean  a sequence with such properties, unless otherwise specified.

 \begin{remark}
 Notice that  sequences with a finite number of elements are excluded from this di\-scus\-sion, but this is a trivial case (both the BM-density and  asymptotic density vanish). Similarly, we are not dealing with sequences  with repetitions (see Remark \ref{remark}). 
\end{remark}

\bigskip
\noindent
The paper is structured as follows. Section 2 contains two equivalent definitions of the BM-density; similarly, Section 3 discusses two equivalent definitions 
 of  the asymptotic density (upper and lower); in Section 4 we prove a general inequality relating the upper asymptotic density and the BM-density; in Section 5 we prove the coincidence  of the asymptotic density and  the BM-density  for a particular class of sequences, to which arithmetic progressions belong. Last, in Section 6 (i) we show how to use the alternative definition \ref{def1} for estimating  the BM-density of some kind of sequences that again generalizes  arithmetic progressions; (ii) we find a nice characterization of  the  BM-density (see Theorem \ref{fourth result} and Corollary \ref{corollary}). 

\bigskip
\noindent
We shall use the standard symbols $\lfloor x\rfloor$ and $\lceil x\rceil$ to mean respectively the greatest integer less than or equal to $x$  ({\it integer part} of $x$) and the least integer greater than or equal to $x$.

  \bigskip
  \noindent
The {\it counting function} of $ \Lambda= \{\lambda_n\}_{n \in \mathbb{N}^*}$ is defined as  
     $$F_\Lambda(t)= \begin{cases}0 & t=0\\ \#\{k  \in \mathbb{N}^*: \lambda_k \leq  t\}& t>0.
     
\end{cases}      $$For $a< b$  we have clearly
$$F_\Lambda(b)-F_ \Lambda(a)= \#\{k  \in \mathbb{N}^*: \lambda_k \in (a, b]  \}.$$

 \section{Various equivalent definitions of the Beurling-Mallia\-vin   density}

 \bigskip 
   \noindent
  The Beurling-Malliavin density $b(\Lambda)$, firstly defined in \cite{BM}, has been studied later in \cite{K2} and in \cite{R} (among others); both the definitions of  \cite{K2} and   \cite{R} are different but equivalent  to the original one. In the very recent paper \cite{GGM}  a further equivalent definition has been given. In the present paper we shall use first the approach of \cite{R}    and later the one of  \cite{GGM}.
  
  \subsection{The definition of  \cite{R} }  In \cite{R} it is proved that an equivalent form of Beurling--Malliavin density is given by the following
  \begin{definition}\label{def1}
 The Beurling Malliavin density  $b(\Lambda)$ of the sequence $\Lambda= (\lambda_n)_{n \in \mathbb{N}^*}$ is the infimum of the numbers $a\geq 0$ for which there exists a sequence of distinct positive integers $N =\big(n _k\big)_{k \in \mathbb{N}^*}$ such that
 $$\sum_{k=1}^\infty \Big|\frac{1}{\lambda_k}- \frac{a}{n _k}\Big|< \infty.$$
\end{definition}   

\bigskip
  \noindent
In \cite{R} it is shown that the set $\mathcal{S}$ of such $a$ is a  right infinite interval  in $\mathbb{R}^+$ (i.e. if $a$ belongs to $\mathcal{S}$ and $b > a$, then $b$ belongs to $\mathcal{S}$ as well). If  $\mathcal{S}$  is empty, then $b(\Lambda)=+\infty.$

\bigskip
\noindent
\begin{remark} \label{8} If $\Lambda $ is a sequence of integers, we have obviously $b(\Lambda)\leq 1$ (take $a=1$ and $N= \Lambda$ in the above definition).

\end{remark}

\subsection{The definition of  \cite{GGM} }
  Let   $\mathfrak{C}$  be the family of all sequences $$\mathcal{I}= \big((a_n, b_n] \big)_{ n \in \mathbb{N}^*}$$
  of intervals in $(0, + \infty)$ such that $a_n < b_n \leq a_{n+1}$ for all $n \in \mathbb{N^*}$ and
  $$\sum_{n = 1}^\infty \Big(\frac{b_n}{a_n}-1\Big)^2= +\infty.$$
  \begin{remark}
  These systems of intervals are referred to as {\rm long} in \cite{P1} and \cite{P2} (see the Introduction).  In \cite{K2} they are called {\rm substantial}. 
  \end{remark}

   \bigskip
  \noindent
  For every $\mathcal{I}= \big((a_n, b_n] \big)_{ n \in \mathbb{N}^*}\in \mathfrak{C} $, let
 \begin{equation}\label{16}
 \ell_{\mathcal{I}}= \liminf_{n \to \infty}\frac{F_\Lambda(b_n) - F_\Lambda(a_n)}{b_n -a_n}.
\end{equation}  

  \bigskip \noindent
 In \cite{GGM} it is proved that

 \begin{proposition}\label{prop9}\sl   
 The following relation holds true
 \begin{equation*}
 b(\Lambda)= \sup\{\ell_{\mathcal{I}},\,\mathcal{I} \in \mathfrak{C} \}.
 \end{equation*}
 \end{proposition}

   \bigskip
  \noindent
 We shall be interested in the   subset  of $\mathfrak{C}$ defined as 
 $$\mathfrak{C}_{>1}= 
   \Big\{ \mathcal{I}= \big((a_n, b_n] \big)_{ n \in \mathbb{N}^*}\in \mathfrak{C} : \limsup_{n \to \infty}\frac{b_n}{a_n  }>1 \Big\}.$$
  Accordingly, we shall denote
  \begin{equation}\label{33} 
  b_{>1}(\Lambda)= \sup\{\ell_{\mathcal{I}},\mathcal{I} \in \mathfrak{C}_{>1} \}.
  \end{equation}  From Proposition \ref{prop9} it follows that
   \begin{equation}b_{>1}(\Lambda)\leq  b(\Lambda).\label{34} 
  \end{equation}

\section{Two equivalent definitions of the asymptotic density}

Let $\Lambda = (\lambda_n)_{n \in \mathbb{N}^*}$ be a sequence of positive numbers. 
\begin{definition}
(i) The {\it lower {\rm (resp.}   upper{\rm )} asymptotic density of $\Lambda$} is  
$$\underline{d}(\Lambda)=\liminf_{n \to \infty} \frac{F_\Lambda(n)}{n}, \qquad  {\rm (resp.}\,\, \overline{d}(\Lambda)=\limsup _{n \to \infty} \frac{F_\Lambda(n)}{n}).$$ 
(ii)  We say that $\Lambda = (\lambda_n)_{n \in \mathbb{N}^*}$ has  {\it asymptotic density} ${d}(\Lambda)$ if $$\underline{d}(\Lambda)=\overline{d} (\Lambda){\rm \big(}={d}(\Lambda){\rm \big)}.$$
\end{definition}

\bigskip
\noindent
It is well known  that, besides the definition,  there is another way for calculating  the upper and lower asymptotic densities of $\Lambda$. Precisely we have
\begin{proposition}\label{prop19}\sl The set of limit points of the sequence $\big(\frac{k}{\lambda_{ k}}\big)_{k \in \mathbb{N}^*}$ coincides with the set of limit points of the sequence $\big(\frac{F_\Lambda(n)}{n}\big)_{n \in \mathbb{N}^*}$. 

\bigskip
\noindent
In particular the following relations hold true
 $$\underline{d}(\Lambda)= \liminf_{k \to \infty }\frac{k}{\lambda_{ k}};   \qquad  \overline{d}(\Lambda)= \limsup_{k \to \infty }\frac{k}{\lambda_{ k}}. $$
 Moreover $d(\Lambda)$ exists if either of the two limits
$$  \lim _{n \to \infty }\frac{F_\Lambda(n)}{n}, \qquad  \lim _{k \to \infty }\frac{k}{\lambda_{ k}}$$
exists; in this case
$$d(\Lambda)=\lim _{n \to \infty }\frac{F_\Lambda(n)}{n}=  \lim _{k \to \infty }\frac{k}{\lambda_{ k}}.$$
\end{proposition}

\bigskip
\noindent
 {\it Proof.}
 In fact, if $\lambda_{k_n} \leq n < \lambda_{ k_n+1}$, we have $F_\Lambda(n)= k_n$, and
\begin{equation}\label{32}
\frac{k_n}{\lambda_{ k_n+1}}<  \frac{F_\Lambda(n)}{n}  \leq \frac{k_n}{\lambda_{ k_n}}.
\end{equation}
 Notice that $\lim_{n \to \infty}k_n = \infty$ and let $\ell$ be a limit point for the sequence $(\frac{F_\Lambda(n)}{n})_{n \in \mathbb{N}^*}$. Then there exists a subsequence $(n_r)_{r \in \mathbb{N}^*}$ converging to $\infty$ such that $\frac{F_\Lambda(n_r)}{n_r}\to \ell$ as $r \to \infty$. Relation \eqref{32} implies that
 $$ \overline{\kappa}:=\limsup_{r \to \infty }\frac{k_{n_r}}{\lambda_{ k_{n_r}}}=\limsup_{r \to \infty }\frac{k_{n_r}}{\lambda_{ k_{n_r}+1}} \leq \ell \leq   \liminf_{r \to \infty }\frac{k_{n_r}}{\lambda_{ k_{n_r}}} =: \underline{\kappa}.$$
In other words every limit point for $\big(\frac{F_\Lambda(n)}{n}\big)_{n \in \mathbb{N}^*}$ is a limit point for $\big(\frac{k}{\lambda_{ k}}\big)_{k \in \mathbb{N}^*}$. 

\bigskip
\noindent
To go the inverse direction, notice that $F_\Lambda(\lambda_k)=k$; thus, for $n_k \leq \lambda_k < n_k+1$,
$$\frac{F_\Lambda(n_k)}{n_k+1}\leq \frac{k}{\lambda_{ k}}\leq \frac{F_\Lambda(n_k+1)}{n_k };$$
now argue   as in the first part of this proof, using the fact that $ \lim_{k \to \infty}n_k= \infty$.

 \B
 
 \section{The general inequality}
 
 \begin{theorem} \label{first result}\sl For every sequence $\Lambda = (\lambda_n)_{n \in \mathbb{N}^*}$, the relation
 \begin{equation}\label{basic inequality}
 \overline d(\Lambda) \leq b(\Lambda)
 \end{equation}
 holds true.
 \end{theorem}
 
 \noindent
 { \it Proof.} Recall that the Beurling--Malliavin Theorem states that, putting$$ \mathcal{E}_{\Lambda}= \{\re^{\pm\ri \lambda_n x}\}, \qquad \mathcal{R}(\Lambda) = \sup \{a \in \mathbb{R}^+: \mathcal{E}_{\Lambda} \hbox{ is complete  in } L^2(0,a)\},$$
 we have
 $$\mathcal{R}(\Lambda) = 2\pi b(\Lambda).$$
 Hence it suffices to prove that, for every $\epsilon > 0$,  $\mathcal{E}_{\Lambda} $ is complete  in  $ L^2\big(0,2\pi  (\overline d(\Lambda)- \epsilon)\big)$.
  Let $f$ be a function in $ L^2\big(0,2\pi  (\overline d(\Lambda)- \epsilon)\big)$ and assume that 
  $$\int_0^{ 2\pi  (\overline d(\Lambda)- \epsilon)}\re ^{\pm\ri \lambda_n x}f(x)\rd x = 0,\qquad n= 1,2, 3, \dots$$
 By the change of variable $y = \frac{x}{\overline d(\Lambda)- \epsilon}- \pi$, 
 the above integral becomes
 $$(\overline d(\Lambda)- \epsilon)\re ^{\pm\ri \lambda_n  (\overline d(\Lambda)- \epsilon) \pi }\int_{-\pi}^\pi \re ^{\pm\ri \lambda_n  (\overline d(\Lambda)- \epsilon)y}f\big( (\overline d(\Lambda)- \epsilon)(y+\pi)\big)\rd y.$$
 In \cite{PW} the following result is proved (see Theorem  XXVIII, p. 84):
     \begin{theorem} \label{PW}\sl Let $0< m_1< m_2 < \cdots $ and let $$\limsup_{n \to \infty} \frac{n}{m_n}>1.$$ Then, if $f \in L^2 (- \pi, \pi)$ and
 $$\int_{-\pi}^\pi \re ^{\pm\ri m_n x}f(x)\rd x = 0,\qquad n= 1,2, 3, \dots\, ,$$
 then $f(x)=0$ for almost every $x\in (-\pi,\pi)$.
 \end{theorem}
 
 \bigskip
 \noindent
 Setting  $m_n =    (\overline d(\Lambda)- \epsilon)\lambda_n$,  we obtain that 
 $$\int_{-\pi}^\pi \re ^{\pm\ri m_n  y}g(y)\rd y=0,$$
 where $g$ is the function  $y \in   (-\pi,\pi)  \mapsto f\big( (\overline d(\Lambda)- \epsilon)(y+\pi)\big)$. 
 
 \medskip
 \noindent
 Observe that
 $$\limsup_{n \to \infty} \frac{n}{m_n} = \limsup_{n \to \infty}\frac{n}{ (\overline d(\Lambda)- \epsilon)\lambda_n}= \frac{\overline d(\Lambda)}{ \overline d(\Lambda)- \epsilon }>1. $$
 Hence Theorem \ref{PW} is in force and  we deduce  that $g$ vanishes almost everywhere on  $(-\pi,\pi)$. This is obviously equivalent to saying that $f$  vanishes almost everywhere on  $2\pi  (\overline d(\Lambda)- \epsilon)$, i.e. that  $\mathcal{E}_{\Lambda} $ is complete  in  $ L^2\big(0,2\pi  (\overline d(\Lambda)- \epsilon)\big)$.

 \B
 
 \begin{remark} Theorem \ref{first result} is stated but not proved in \cite{P}. In any case, \cite{P} doesn't formulate it in terms of densities; see slide 3 in \cite{P}.
 \end{remark}

 \noindent
 The following example shows that the inequality \eqref{basic inequality} may be strict.
 
 \begin{example} Let $p$  be an integer with $ 1\leq p \leq 9$ and 
denote by  $\Lambda = (\lambda_n)_{n \in \mathbb{N}^*}$ the - strictly increasing - sequence formed by the positive integers having first digit equal to $p$.  It is well known (see for instance \cite{FG}) that $\overline d(\Lambda)= \frac{10}{9(p+1)}<1$, and now we prove that $b(\Lambda)=1$. Recall that  the BM-density of any set of integers is $\leq 1$. By formula \eqref{33} and relation \eqref{34}, it suffices to prove that there exists a family of intervals
    $\mathcal{I}= \big((a_n, b_n] \big)_{ n \in \mathbb{N}^*} $ such that  $\limsup_{n \to \infty}\frac{b_n}{a_n  }>1  $ and  $$
 \ell_{\mathcal{I}}= \liminf_{n \to \infty}\frac{F_\Lambda(b_n) - F_\Lambda(a_n)}{b_n -a_n} \geq 1.$$
 Take
 $$a_n = p 10^n, \qquad b_n =(p+1)10^n-1.$$
 We have $b_n- a_n=10^n-1$ and 
 $$F_\Lambda(b_n) =\sum_{k=1}^n \sum_{j= p 10^k}^{(p+1)10^k-1}1= \sum_{k=1}^n \{(p+1)10^k- p 10^k\}= \sum_{k=1}^n 10^k= \frac{10}{9}(10^n-1); $$
 $$F_\Lambda(a_n) = F_\Lambda(b_{n-1})+ 1= \frac{10}{9}(10^{n-1}-1)+ 1=  \frac{1}{9}(10^n-1).$$
 Thus
 $$\ell_{\mathcal{I}}= \liminf_{n \to \infty}\frac{F_\Lambda(b_n) - F_\Lambda(a_n)}{b_n -a_n}= \liminf_{n \to \infty}\frac{ \frac{10}{9}(10^n-1)- \frac{1}{9}(10^n-1)}{10^n-1}=1 .$$
 \end{example}
 
 \begin{remark} It is easy to check that for $a_n=(p+1)10^n-1 $ and $b_n = a_{n+1}= (p+1)10^{n+1}-1 $ we have
 $$\ell_{\mathcal{I}}=\overline d(\Lambda)= \frac{10}{9(p+1)}.$$
 \end{remark}
  
 \section{A general result concerning the equality}  
 The aim of this Section is to prove the following result, in which a class of sequences $\Lambda$ is identified for which $d(\Lambda)= b(\Lambda)$.
 \begin{theorem}\label{second result}\sl Assume that $\Lambda = (\lambda_n)_{n \in \mathbb{N}^*}$ is a sequence with the following property: there exists an integer $p \geq 1$ and non-negative constants $a_0, a_1, \dots, a_{p-1}$ such that
 $$\lim_{k\to \infty}\{\lambda_{kp +j+1}-\lambda_{kp +j}\} = a_j, \quad \, \forall j = 0,1, 2, \dots p-1.$$
 Assume in addition that $\sum_{j = 0}^{p-1} a_j>0$. Then 
 $$ d(\Lambda)= b(\Lambda)= \frac{p}{\sum_{j = 0}^{p-1} a_j}.$$
 \end{theorem}
 
 \bigskip
 \noindent
 For the proof of this result, we need some lemmas.
 \begin{lemma}\label{lemma1}\sl Assume that  $\Lambda = (\lambda_n)_{n \in \mathbb{N}^*}$ verifies the conditions of Theorem \ref{second result}. Then
 $$\lim_{n \to \infty} \frac{\lambda_n}{n}= \frac{\sum_{j = 0}^{p-1} a_j}{p}.$$
 \end{lemma}
 
 \noindent
 { \it Proof.} We shall prove that, for every sequence $(\lambda_{kp +j})_{k \in \mathbb{N}^*} $, $j=0, 1, \dots, p-1$, we have
 $$\lim_{k\to \infty}\frac{\lambda_{kp +j}}{kp +j}= \frac {\sum_{j = 0}^{p-1} a_j}{p}.$$
 First, by Cesaro Theorem, 
 $$\lim_{k\to \infty}\frac{\lambda_{kp +j}}{kp +j}= \lim_{k\to \infty} \frac{\lambda_{(k+1) p +j}-\lambda_{kp +j}}{ (k+1) p +j- (kp +j)}= \lim_{k\to \infty} \frac{\lambda_{(k+1) p +j}-\lambda_{kp +j}}{ p}.$$
 Now
 \begin{align*}&
 \lim_{k\to \infty}\{\lambda_{(k+1) p +j}-\lambda_{kp +j}\}= \sum_{r=j}^{p+j-1}\lim_{k\to \infty}\{\lambda_{kp +r+1}-\lambda_{kp +r}\}\\ &=  \sum_{r=j}^{p-1}\lim_{k\to \infty}\{\lambda_{kp +r+1}-\lambda_{kp +r}\} + \sum_{r=p}^{p+j-1}\lim_{k\to \infty}\{\lambda_{kp +r+1}-\lambda_{kp +r}\}\\ &= 
\sum_{r=j}^{p-1}  a_r+ \sum_{s=0}^{j-1}\lim_{k\to \infty}\{\lambda_{(k+1)p +s+1}-\lambda_{(k+1)p +s}\}\\& =  \sum_{r=j}^{p-1}  a_r+ \sum_{s=0}^{j-1}\lim_{k\to \infty}\{\lambda_{kp +s+1}-\lambda_{kp +s}\}\\&=  \sum_{r=j}^{p-1}  a_r+\sum_{s=0}^{j-1}a_s=\sum_{r = 0}^{p-1} a_r.\end{align*}

 \B 
 
 \begin{remark} \label{remark1}Lemma \ref{lemma1} implies that any sequence $\Lambda= (\lambda_n)_{n \in \mathbb{N}^*}$ satisfying the assumptions of Theorem \ref{second result} has the form $$\lambda_n = \frac{\sum_{j = 0}^{p-1} a_j}{p} n + \psi_n,$$
 where 
 \begin{equation}\label{7}
 \lim_{n \to \infty } \frac{\psi_n}{n}=0.
 \end{equation}
Besides \eqref{7}, a necessary condition for $\Lambda$ to satisfy the assumptions of Theorem \ref{second result} is that
$$\lim_{n \to \infty}\psi_{n+p}-\psi_n=0.$$
This is easily seen:
$$\psi_{n+p}-\psi_n = \lambda_{n+p} -\frac{\sum_{j = 0}^{p-1} a_j}{p} (n+p) - \lambda_n +   \frac{\sum_{j = 0}^{p-1} a_j}{p} n =  \lambda_{n+p}- \lambda_n-\sum_{j = 0}^{p-1} a_j=0 $$
since, for $n = kp + j$, $j =0, \dots, p-1$, we have
 $$ \lim_{k \to \infty}\{\lambda_{kp + j+p} - \lambda_{kp + j}\}= \sum_{j = 0}^{p-1} a_j, $$
 as  has been proven in Lemma \ref{lemma1}.\end{remark}
 
 \begin{lemma}\label{28}\sl Let $(x_n)_{n \in \mathbb{N}^*}$ and $(y_n)_{n \in \mathbb{N}^*}$ be two sequences of   numbers, with 
 \begin{equation}\label{29}
  \lim_{n \to \infty}(y_n - x_n)= + \infty.
 \end{equation}
  Then $$\lim_{n \to \infty} \frac{\lfloor y_n\rfloor-\lfloor x_n\rfloor}{y_n - x_n}=1.$$
\end{lemma}

\bigskip
\noindent{\it Proof}.  Let $\{x\}= x- \lfloor x \rfloor$ be the fractional part of the number $x$ and   write
$$\frac{\lfloor y_n\rfloor-\lfloor x_n\rfloor}{y_n - x_n}= \frac{1}{1 + \frac{\{y_n\}-\{x_n\}}{\lfloor y_n\rfloor-\lfloor x_n\rfloor}}.$$
Observe that  $ |\{y_n\}-\{x_n\}|< 1$ (from $0 \leq \{x\}< 1$) and that $\lim_{n \to \infty}( \lfloor y_n\rfloor-\lfloor x_n\rfloor)= + \infty$ (since $\lfloor y_n\rfloor-\lfloor x_n\rfloor\geq  y_n - x_n-1$). The statement follows.
 
\B

\bigskip
\noindent
 Denote by $\overline{\lambda}$ the interpolated linear function relative to $ \Lambda$, i.e. the function defined by
$$\overline{\lambda}(x) = \big(\lambda_{k +1}-\lambda_k\big)(x-k)+ \lambda_k, \qquad k \leq x < k+1, \qquad k=1, 2 \dots$$
or equivalently
$$\overline{\lambda}(x) =\big(\lambda_{\lceil x \rceil}-\lambda_{ \lfloor x \rfloor}\big)(x-\lfloor x \rfloor)+ \lambda_{\lfloor x \rfloor}.$$

\bigskip
\noindent
From \cite{BGT} we recall the
\begin{definition} The function $f: \mathbb{R}^+\to  \mathbb{R}^+$ is {\it regularly varying with exponent $\rho>0$} if, for every $t>0$, we have
$$\lim_{x \to +\infty}\frac{f(tx)}{f(x)}= t^\rho.$$

\end{definition}

\begin{lemma}\label{lemma2}\sl Assume that  $\Lambda = (\lambda_n)_{n \in \mathbb{N}^*}$ is such that
\begin{equation}\label{13}
\lim_{n \to \infty} \frac{\lambda_n}{n}=\ell >0.
\end{equation}
Then
\begin{itemize} \item[(i)]$\lim_{x \to +\infty} \frac{\overline{\lambda}( x  )}{x} = \ell;$

\item[(ii)] $\overline{\lambda}$ is regularly varying with exponent 1;
\item[(iii)] for every $\beta > 1$ we have
$$\lim_{ x\to +\infty \atop
  \beta x< y \to +\infty }\frac{\overline\lambda(y)-\overline\lambda(x)}{y-x}= \ell .$$

\end{itemize} \end{lemma}

 \noindent
{ \it Proof.}  First,  from the assumption \eqref{13} and the relations
$$  \frac{\lfloor x \rfloor}{\lceil x \rceil}\cdot \frac{\overline{\lambda}(\lfloor x \rfloor)}{ \lfloor x \rfloor}= \frac{\overline{\lambda}(\lfloor x \rfloor)}{ \lceil x \rceil}\leq \frac{\overline{\lambda}(\lfloor x \rfloor)}{  x  }\leq \frac{\overline{\lambda}(\lceil x \rceil)}{ \lfloor x \rfloor } $$
and 
$$ \frac{\overline{\lambda}(\lceil x \rceil)}{ \lceil x \rceil}\leq \frac{\overline{\lambda}(\lceil x \rceil)}{  x  }\leq \frac{\overline{\lambda}(\lceil x \rceil)}{ \lfloor x \rfloor }=\frac{\overline{\lambda}(\lceil x \rceil)}{ \lceil x \rceil}\cdot\frac{\lceil x \rceil}{\lfloor x \rfloor},$$ (where we have used also the fact that $\overline{\lambda}$ is non-decreasing) and by the  relation $\lim_{x \to + \infty}\frac{\lceil x \rceil}{\lfloor x \rfloor}=1$, we obtain
\begin{equation}\label{15}
\lim_{x \to +\infty} \frac{\overline{\lambda}( \lfloor x \rfloor   )}{x} =\lim_{x \to +\infty} \frac{\overline{\lambda}( \lceil x \rceil )}{x} =\ell,
\end{equation}    
 and (i) is obvious, since
 $$\frac{\overline{\lambda}( \lfloor x \rfloor   )}{x} \leq \frac{\overline{\lambda}( x  )}{x}\leq\frac{\overline{\lambda}( \lceil x \rceil )}{x}.$$
 For (ii), we have to prove that  for every $t >0$,
 $$\lim_{x \to \infty}  \frac{\overline{\lambda}(t x )}{\overline{\lambda}( x ) }=t,$$ 
 and this can be obtained  by writing
 $$ \frac{\overline{\lambda}(t x )}{\overline{\lambda}( x ) } = t \cdot \frac{\overline{\lambda}( t x )}{t x}\cdot \frac{x}{\overline\lambda( x )}$$
 and using point (i).
 
 \medskip
 \noindent
 Now we prove point (iii). Write
 $$\frac{\overline\lambda(y)-\overline\lambda(x)}{y-x}=  \frac{\big(\lambda_{\lceil y \rceil}-\lambda_{ \lfloor y \rfloor}\big)(y-\lfloor y \rfloor)}{y-x}+  \frac{\big(\lambda_{\lceil x \rceil}-\lambda_{ \lfloor x \rfloor}\big)(x-\lfloor x \rfloor)}{y-x}+\frac{\lambda_{\lfloor y \rfloor}-\lambda_{\lfloor x \rfloor}}{y-x}.$$
 Notice that
 $$\limsup_{x \to +\infty}\big(\lambda_{\lceil x \rceil}-\lambda_{ \lfloor x \rfloor}\big)(x-\lfloor x \rfloor)< +\infty.$$
 This is because $x-\lfloor x \rfloor= \{x\}$ (fractional part of $x$) is between $0$ and $1$; furthermore
 $$\limsup_{x \to +\infty}\big(\lambda_{\lceil x \rceil}-\lambda_{ \lfloor x \rfloor}\big)= \max_{0 \leq j \leq p-1}a_j.$$
 Hence the first two summands above converge to 0, since $y-x > x(\beta -1)\to + \infty$ as $x\to + \infty$.
 
 \noindent
 Concerning the last summand, we write it as
 $$ \frac{\lambda_{ \lfloor y \rfloor}}{y}+ \Big(\frac{\lambda_{ \lfloor y \rfloor}}{y}-\frac{\lambda_{ \lfloor x \rfloor}}{x}\Big)\cdot \frac{x}{y-x}$$
which allows to obtain  the claimed result by \eqref{15}  and the fact that $$0 \leq \frac{x}{y-x} \leq \frac{1}{\beta -1}.$$
 \B
 
 \bigskip
\noindent
 {\it Proof of Theorem \ref{second result}}.
 
 \bigskip
\noindent
Recall that the   asymptotic density  $d(\Lambda)$ of the sequence $\Lambda$  is  the $\lim_{n \to \infty} \frac{n}{\lambda_n} $,
if this limit exists (see Proposition \ref{prop19}). Hence by Lemma \ref{lemma1} we have $d(\Lambda) = \frac{p}{\sum_{j = 0}^{p-1} a_j} $, and now we shall prove  that
\begin{equation}\label{12}
\frac{p}{\sum_{j = 0}^{p-1} a_j} \leq b_{>1}(\Lambda). 
\end{equation}
Recalling \eqref{33} and  the formula preceding \eqref{33}, let $(a_n)_{n \in \mathbb{N}^*}$ and $(b_n)_{n \in \mathbb{N}^*}$ be such that $\lim_{ n\to \infty} {a_n }= +\infty$ and $$\limsup_{ n\to \infty} \frac{b_n}{a_n }\geq \liminf_{ n\to \infty} \frac{b_n}{a_n }= \alpha >1;$$  we are interested in calculating
$\frac{F_\Lambda(b_n)-F_\Lambda(a_n) }{b_n - a_n}$ (see relation \eqref{16}).

\bigskip
\noindent
For every $x \in \mathbb{R}^+$ with $  \lambda_k \leq x <\lambda_{k+1}$ we have $F_\Lambda(x)=k $; since $k \leq \big(\overline{\lambda}\big)^{-1}(x)< k+1$, we can write $F_\Lambda(x)= \lfloor \big(\overline{\lambda}\big)^{-1}(x)\rfloor $.

\bigskip
\noindent
  Thus
$$\frac{F_\Lambda(b_n)-F_\Lambda(a_n) }{b_n - a_n}= \frac{\lfloor \big(\overline{\lambda}\big)^{-1}(b_n)\rfloor-\lfloor \big(\overline{\lambda}\big)^{-1}(a_n)\rfloor}{b_n - a_n}= \frac{\lfloor y_n \rfloor-\lfloor x_n \rfloor }{\overline{\lambda}(y_n) -\overline{\lambda}(x_n)  }=  \frac{y_n - x_n  }{\overline{\lambda}(y_n) -\overline{\lambda}(x_n)}\cdot\frac{\lfloor y_n\rfloor-\lfloor x_n\rfloor}{y_n - x_n} ,$$
where
$$x_n =  \big(\overline{\lambda}\big)^{-1}(a_n), \qquad y_n =  \big(\overline{\lambda}\big)^{-1}(b_n).$$
Notice that $x_n \to + \infty$. Then \eqref{12}  follows immediately from \eqref{34},  Lemma \ref{lemma2} (iii) and Lemma \ref{28} (see \eqref{29}) if we prove that there exists $\beta > 1$ such that $y_n > \beta x_n$ for every sufficiently large $n$.

\bigskip
\noindent
Let $\epsilon>0$ be fixed with $\epsilon < \alpha -1$. There exists $n_0$ such that, for $n>n_0$, $(\alpha - \epsilon)a_n\leq b_n$, which yields $$\big(\overline{\lambda}\big)^{-1}\big( (\alpha - \epsilon)a_n\big)\leq \big(\overline{\lambda}\big)^{-1}(b_n).$$
By Lemma \ref{lemma2} (ii) and Theorem 1.5.12 in \cite{BGT}, $\big(\overline{\lambda}\big)^{-1}$ is regularly varying with exponent 1,  hence
$$\frac{\big(\overline{\lambda}\big)^{-1}( (\alpha - \epsilon)a_n)}{\big(\overline{\lambda}\big)^{-1}(a_n)}\to \alpha - \epsilon , \qquad n \to \infty.$$
It follows that for every $\delta > 0$   with $\delta <  \alpha - \epsilon -1 $, we have
$$ \big(\overline{\lambda}\big)^{-1}(a_n)\Big( \alpha - \epsilon  - \delta \Big) \leq \big(\overline{\lambda}\big)^{-1}\big( (\alpha - \epsilon)a_n\big)\leq \big(\overline{\lambda}\big)^{-1}(b_n)$$
or $y_n \geq \beta x_n $ with $\beta = \alpha - \epsilon - \delta >1$.

\bigskip
\noindent
Now an application of Theorem \ref{first result} concludes the proof.

\B
 
\begin{example}  \label{example} Arithmetic progressions satisfy  the assumption of Theorem \ref{second result}. A less trivial example is
$$\lambda_n = \frac{3}{2}n + \frac{3 + (-1)^n}{4};$$
notice that the values of this sequence are all integers. Hence $$b(\Lambda) = d(\Lambda) = \frac{2}{3}.$$
\end{example}

  \section{Some results that follow from the definition of \cite{R}}  
  In this Section we show how to use Definition \ref{def1} to get information about $b(\Lambda)$ and its connection with $d(\Lambda)$.
\begin{proposition}\label{third result}\sl
Assume that there exist  $\ell >0$ and a  sequence  
 $(\psi_n) _{n \in \mathbb{N}^*}  $  such that $$\lambda_n= \ell n + \psi_n  .$$ If either
 
 \noindent
(i)  $ (\psi_n) _{n \in \mathbb{N}^*} $ is non-decreasing (ultimately)
 
 \medskip
 \noindent
 or
 
 \medskip
 \noindent
 (ii)   $  \kappa :=\liminf_{n \to \infty} \Big|\ell+\frac{ \psi_n}{n }\Big|>   0 $   and $\sum_n\frac{|\psi_n|}{n^2}< + \infty,$  
 
 \medskip
 \noindent
 then
\begin{equation}\label{14}b(\Lambda)\leq  \frac{1}{\ell}.
\end{equation}
Moreover, in both cases, we have $\alpha:= \liminf_{n \to \infty}\frac{ \psi_n}{n } \geq 0$. In particular, if   $\alpha=0,$ then
$$b(\Lambda)=  \overline d(\Lambda) = \frac{1}{\ell }.$$
\end{proposition}

 \noindent
 {\it Proof.}  We   prove that,    in both cases (i) and (ii),    the number $\frac{1}{\ell}$ verifies the property stated in   Definition \ref{def1}.

 \bigskip
 \noindent
 (i) Let $n_k = k + \big\lfloor \frac{ \psi_k}{\ell}\big\rfloor$. This sequence of integers is made of distinct integers since 
 $$n_{k+1}-n_k= \Big\lfloor \frac{ \psi_{k+1}}{\ell}\Big\rfloor-\Big\lfloor \frac{ \psi_k}{\ell}\Big\rfloor+ 1 \geq 1.$$
 Moreover
 \begin{align*}&
 \sum_k \Big|\frac{1}{\lambda_k}- \frac{\frac{1}{\ell}}{n_ k}\Big| =\frac{1}{\ell}\sum_k \Big| \frac{1}{k + \frac{\psi_k}{\ell}}- \frac{1}{ k + \lfloor\frac{\psi_k}{\ell}\rfloor }\Big|= \frac{1}{\ell}\sum_k  \frac{|\frac{\psi_k}{\ell}-\lfloor\frac{\psi_k}{\ell}\rfloor|}{(k + \frac{\psi_k}{\ell})( k + \lfloor\frac{\psi_k}{\ell}\rfloor)}\\&\leq\frac{1}{\ell}\sum_k \frac{1}{k^2(1 + \frac{\psi_k}{\ell k})( 1 + \frac{1}{k}\lfloor\frac{\psi_k}{\ell}\rfloor)}.
 \end{align*}
Since $(\psi_n)_{n \in \mathbb{N}^*} $ is non-decreasing, there exists an integer $k_0$ and a constant $\alpha$ such that for $k \geq k_0$,   $$ 1 + \frac{\alpha}{ k}>0, \qquad 1 + \frac{\psi_k}{\ell k}\geq 1 + \frac{\alpha}{ k}, \qquad  1 + \frac{1}{k}\Big\lfloor\frac{\psi_k}{\ell}\Big\rfloor\geq 1 + \frac{\alpha}{ k}.$$ 
Thus the series  $\sum_k \frac{1}{k^2(1 + \frac{\psi_k}{\ell k})( 1 + \frac{1}{k}\lfloor\frac{\psi_k}{\ell}\rfloor)}$ is majorized by $$const. +  \sum_k \frac{1}{k^2(1 + \frac{\alpha}{  k})^2}, $$ 
which behaves like $\sum_k \frac{1}{k^2}< + \infty.$

  \bigskip
 \bigskip
 \noindent  (ii) Take $n_k =k$. We have $$\sum_k \Big|\frac{1}{\lambda_k}- \frac{\frac{1}{\ell}}{ k}\Big| = \frac{1}{\ell^2}\sum_k  \frac{|\psi_k  |}{k^2|1+ \frac{\psi_k}{\ell k}| } $$
which, by the assumption $\kappa >0$,  behaves like $\sum_k  \frac{|\psi_k  |}{k^2 } < + \infty.$

 \bigskip
 \noindent 
Concerning  $\overline d(\Lambda) $, we have
$$\overline d(\Lambda) = \limsup_{n \to \infty} \frac{n}{\lambda_n}= \limsup_{n \to \infty} \frac{1}{\ell + \frac{\psi_n}{n}}= \frac{1}{\ell +\limsup_{n \to \infty}  \frac{\psi_n}{n} } = \frac{1}{\ell+\alpha}.$$

\bigskip
\noindent
The inequality $\alpha \geq 0$ is obvious in case (i);   in case (ii) it follows from the basic inequality \eqref{basic inequality}.
If $\alpha=0$,  the above discussion implies that
$$b(\Lambda)\leq \frac{1}{\ell}=\overline d(\Lambda),$$
and we get  the equality once more by \eqref{basic inequality}.

\B

\begin{example}   Any arithmetic progression $\lambda_n  = cn + a$  verifies both the assumptions (i) and (ii) of Proposition \ref{third result}; in particular the set of multiples of the integer $p$ has Beurling-Malliavin density equal to $\frac{1}{p}$.  

\noindent
The sequence of Example \ref{example} verifies (ii) of Proposition \ref{third result}. 
\end{example}

\begin{remark} There are sequences which satisfy   the assumptions of Proposition \ref{third result} but not those of Theorem \ref{second result}. For instance take  
$$\lambda_n  = \ell n + \sin n.$$
By Proposition \ref{third result} (ii) 
$$b(\Lambda)= d(\Lambda) = \frac{1}{\ell }$$
($d(\Lambda)$ exists since $\lim_{n \to \infty}\frac{\sin n}{n}=0$). Theorem \ref{second result} doesn't hold since, for every $p$, the sequence $\big(\sin (n + p)- \sin n\big)_{n \in \mathbb{N}}$ doesn't have any limit   (see Remark \ref{remark1}).

\end{remark}  
\begin{proposition}\label{prop8} \sl Assume that $\liminf_{n \to \infty}(\lambda_{n+1}-\lambda_n) =: \ell>0$. Then
$$b(\Lambda) \leq \frac{1}{\ell}.$$

\end{proposition}
 
  \bigskip
 \noindent
  {\it Proof.}  Recall that $\Lambda$ is strictly increasing, hence $\ell \geq 0$. Now let $\ell>0$.  It suffices to prove that $$b(\Lambda) \leq \frac{1}{\ell - \epsilon }$$
for every $\epsilon$ with $0<\epsilon < \ell$.  

\medskip
\noindent
Take $n_k= \lfloor \frac{\lambda_k}{ \ell-\epsilon}\rfloor$. Then
  $$n_{k + 1} -n_k =  \Big\lfloor \frac{\lambda_{k+1}}{ \ell-\epsilon}\Big\rfloor-\Big\lfloor \frac{\lambda_k}{ \ell-\epsilon}\Big\rfloor\geq \frac{\lambda_{k+1}}{ \ell-\epsilon}-1- \frac{\lambda_k}{ \ell-\epsilon}>0$$
ultimately, since $\lambda_{k+1}-\lambda_k> \ell- \epsilon$ ultimately. Hence the sequence $(n_k)_{k \in \mathbb{N}^*}$ is made of distinct integers (ultimately). 

\medskip
\noindent
 Moreover\begin{align*}&
\sum_k \Big|\frac{1}{\lambda_k}- \frac{\frac{1}{\ell - \epsilon}}{ n_k}\Big| =  \sum_k \Big|\frac{1}{\lambda_k}- \frac{\frac{1}{\ell - \epsilon}}{ \lfloor \frac{\lambda_k}{\ell - \epsilon} \rfloor}\Big| =  \sum_k \frac{|\lfloor \frac{\lambda_k}{\ell - \epsilon} \rfloor- \frac{\lambda_k}{\ell - \epsilon}|}{\lambda_k   \lfloor \frac{\lambda_k}{\ell - \epsilon}  \rfloor} \leq \sum_k \frac{1}{\lambda_k\lfloor  \frac{\lambda_k}{\ell - \epsilon}  \rfloor},
\end{align*} 
  and the last series has the same behaviour as $\sum_k  \frac{1}{\lambda^2_k} $ (since $\lambda_k \to \infty$). It remains to prove that 
 $$\sum_k  \frac{1}{\lambda^2_k}  <+ \infty.$$
  Putting $\psi_n =\lambda_{n+1}-\lambda_n  $, it is easy to see by induction that
  $$\lambda_ n= \lambda_ 1 + \sum_{k= 1}^{n-1}\psi_k.$$
 As stated above, there exists $n_0$ such that $\psi_n> \ell - \epsilon$ for $n > n_0$.
  Hence, for $n > n_0+1$, we have
  $$\lambda_ n= \lambda_ 1 + \sum_{k= 1}^{n_0}\psi_k + \sum_{k= n_0+1}^{n-1}\psi_k >  \lambda_ 1 + \sum_{k= 1}^{n_0}\psi_k +(\ell - \epsilon) (n-1-n_0) = (\ell - \epsilon) n + const., $$
 which implies that
  $\sum_k  \frac{1}{\lambda^2_k}$ is  $ \leq const\cdot\sum_k  \frac{1}{k^2 }<+ \infty.$

  \B

\bigskip
\noindent
The preceding Proposition and   Theorem \ref{first result} yield immediately the following result, which is a particular case of (actually it has been a motivation for)  Theorem \ref{second result}:
\begin{corollary} \sl If the limit
$$ \ell:=\lim_{n \to \infty}(\lambda_{n+1}- \lambda_n)$$
exists (as an extended number)   and is strictly positive, then
$$d(\Lambda)= b(\Lambda) = \frac{1}{\ell},$$
where we adopt the convention $\frac{1}{0}= + \infty$. 

\end{corollary}

\begin{remark} \label{remark} For the sequence $\Lambda =(\lambda_n)_{n \in \mathbb{N}^*}$ with $\lambda_n = \log n$ we have $\ell =0$ and as a consequence $b(\Lambda)=+\infty$. We conjecture that the same happens for $M =(\mu_n)_{n \in \mathbb{N}^*}$ with $\mu_n = \lfloor\log n\rfloor$, as well as for sequences with \lq\lq too many\rq\rq\ repetitions; but we have no proof at present for this.

\end{remark}

\bigskip
\noindent
 \begin{proposition}\label{prop18} \sl Assume that  $\Lambda =(\lambda_n)_{n \in \mathbb{N}^*}$ is a sequence   such that $$\sum_k  \frac{1}{\lambda_k}  <+ \infty.$$ Then $d(\Lambda) = b(\Lambda)=0$.
 
 \end{proposition} 
  \begin{remark}  This Proposition was proved in \cite{S}, but only in the part concerning $d(\Lambda)$.
 \end{remark}
 
 \bigskip
 \noindent
 {\it Proof.} Recall Olivier's Theorem (also known as  as Abel's or   Pringsheim's Theorem), proved in \cite{O}:
 \begin{theorem} \sl Let $(a_n)_{n\in \mathbb{N}^*}$ be a non--increasing sequence of positive numbers
such that the corresponding series $\sum_{n= 1}^\infty a_n$ is convergent. Then $\lim_{n \to \infty} na_n = 0.$
\end{theorem}
 
 \bigskip
 \noindent
  Applying it with $ a_n =\frac{1}{\lambda_n}$ we have  $$d(\Lambda) =\lim_{n \to \infty} \frac{n}{\lambda_n }=0.$$
 Moreover, by Definition \ref{def1}, $ b(\Lambda)     = 0  $.  
 
 \B

 \begin{example} 
 This result applies for instance to the set of powers $(k^p)_{k \in \mathbb{N}^*}$ ($p \in \mathbb{N}^*$, $p \geq 2 $). More generally, it applies to any  increasing sequence  of the form $(L(k)k^p)_{k \in \mathbb{N}^*}$, where $(L_k) _{k \in \mathbb{N}^*}$ is a slowly varying sequence: recall that a slowly varying sequence $(L_k) _{k \in \mathbb{N}^*}$ has the property that, for any $\epsilon > 0$, ultimately we have  $k^{- \epsilon }<L_k$ (this is an easy consequence of Theorem 1.3.1 in \cite{BGT}), which implies that $\sum_k  \frac{1}{\lambda_k}  < \sum_k  \frac{1}{ k^{p-\epsilon}}< +\infty $ if $\epsilon< p-1 $.
 \end{example}
  
 \bigskip
\noindent
Last, here below is the announced characterization of the BM-density:
 \begin{theorem} \label{fourth result}\sl Let  $\Lambda =(\lambda_n)_{n \in \mathbb{N}^*}$ be a sequence with $b(\Lambda)\in (0, + \infty)$; assume that  there exist   a strictly positive number $\ell$ and  another sequence $M = (\mu_n)_{n \in \mathbb{N}^*}$  with   $b(M)\in (0, + \infty)$ such that 

$$\sum_k\Big|\frac{1}{\lambda_k\ } - \frac{1}{\ell\mu_k}\Big|< +\infty.$$
Then $\ell$ is unique; furthermore
\begin{equation}\label{35}
b(M) = \ell \,b(\Lambda).
\end{equation}

\end{theorem}

  \bigskip
 \noindent
 {\it Proof.} Denote by $$\mathcal A=\Big\{\ell>0:\sum_k\Big|\frac{1}{\lambda_k\ } - \frac{1}{\ell\mu_k}\Big|< +\infty \Big\}$$
 and $$\sigma= \inf\mathcal A, \qquad \tau= \sup\mathcal A.$$
 Let $\epsilon > 0$ and  $N= (n_k)_{k \in \mathbb{N}^*}$ a sequence of distinct integers such that
 \begin{equation} 
 \sum_k \Big|\frac{1}{\mu_k}- \frac{b(M)+ \epsilon}{n_k}\Big|< +\infty.
 \end{equation} 
  Then, for every $\ell \in \mathcal A $,
  \begin{align}&\nonumber
  \sum_k \Big|\frac{1}{\lambda_k}- \frac{b(M)+ \epsilon}{\ell n_k}\Big|\leq   \sum_k\Big|\frac{1}{\lambda_k\ } - \frac{1}{\ell\mu_k}\Big|+  \frac{1}{\ell}\sum_k\Big|\frac{1}{\mu_k}- \frac{b(M)+ \epsilon}{n_k}\Big|< + \infty.
  \end{align}
This yields
  $$b(\Lambda) \leq \frac{ b(M)+ \epsilon}{\ell}$$
 and  by optimizing in $\ell$ and $\epsilon$ we conclude that
 \begin{equation}\label{9}
 b(\Lambda) \leq \frac{ b(M) }{\tau}.
 \end{equation}
 Relation \eqref{9}, together with the assumption $b(\Lambda)>0$, implies that $\tau< \infty$.
 
 \bigskip
\noindent
Now interchange the roles of $\Lambda$ and $M$; we obtain that
$$b(M) \leq \ell  \big( b(\Lambda)+ \epsilon\big), $$
 and optimizing
 \begin{equation}\label{10}
 b(M) \leq \sigma   b(\Lambda).
 \end{equation}
 From relation \eqref{10} and the assumption $b(M)>0$ we deduce that $\sigma>0$; and now, putting together \eqref{9} and \eqref{10}, we get
$$b(\Lambda) \leq \frac{ b(M) }{\tau}\leq \frac{\sigma }{\tau}b(\Lambda);$$
this proves the unicity of $\ell$ and the relation \eqref{35}.
 
 \B
 
 \bigskip
\noindent
 We put in evidence the following particular case (see also Proposition \ref{prop18} for the case $\ell =0$):

\begin{corollary}\label{corollary}\sl Let  $\Lambda =(\lambda_n)_{n \in \mathbb{N}^*}$ be a sequence with $b(\Lambda)< +\infty$; if there exists  a positive number $\ell$ such that
$$\sum_k \Big|\frac{1}{\lambda_k}- \frac{\ell}{ k}\Big|< \infty.$$
then $\ell= b(\Lambda)$.
\end{corollary}

 \end{document}